\documentclass{article}

\usepackage{PRIMEarxiv}

\usepackage[utf8]{inputenc} % allow utf-8 input
\usepackage[T1]{fontenc}    % use 8-bit T1 fonts
\usepackage{hyperref}       % hyperlinks
\usepackage{url}            % simple URL typesetting
\usepackage{booktabs}       % professional-quality tables
\usepackage{amsfonts}       % blackboard math symbols
\usepackage{nicefrac}       % compact symbols for 1/2, etc.
\usepackage{microtype}      % microtypography
\usepackage{lipsum}
\usepackage{fancyhdr}       % header
\usepackage{graphicx}       % graphics
\graphicspath{{media/}}     % organize your images and other figures under media/ folder

\usepackage{amsmath}

\newcommand{\R}{{\mathbb R}}

\newcommand{\dSz}{\mathrm{d} S_z}

\newcommand{\dx}{\mathrm{d}x}

\newcommand{\dr}{\mathrm{d}r}
\newcommand{\norm}[1]{\left\|{#1}\right\|}
\newcommand{\inner}[2]{\langle {#1}, {#2} \rangle}

\newtheorem{theorem}{Theorem}[section]
\newtheorem{lemma}[theorem]{Lemma}
\newtheorem{remark}[theorem]{Remark}
\newtheorem{proposition}[theorem]{Proposition}
\newtheorem{definition}[theorem]{Definition}
\newtheorem{example}{Example}[section]

\title{Optimal Transport and Wasserstein Barycenter for Radially Contoured Distributions
}

\author{
  Keyu Chen \\
  School of Mathematical Sciences, Fudan University\\
  Shanghai 200433, P. R. China\\
  \texttt{chenky21@m.fudan.edu.cn} \\
   \And
  Yunxin Zhang \\
  School of Mathematical Sciences, Fudan University \\
  Shanghai Key Laboratory for Contemporary Applied Mathematics \\
  Shanghai 200433, P. R. China\\
  \texttt{xyz@fudan.edu.cn} \\
}

\begin{document}
\maketitle

\begin{abstract}
	The optimal transport and Wasserstein barycenter of Gaussian distributions have been solved. In literature, the closed form formulas of the Monge map, the Wasserstein distance and the Wasserstein barycenter have been given. Moreover, when Gaussian distributions extend more generally to elliptically contoured distributions, similar results also hold true. In this case, Gaussian distributions are regarded as elliptically contoured distribution with generator function $e^{-x/2}$. However, there are few results about optimal transport for elliptically contoured distributions with different generator functions. In this paper, we degenerate elliptically contoured distributions to radially contoured distributions and study their optimal transport and prove their Wasserstein barycenter is still radially contoured. For general elliptically contoured distributions, we give two numerical counterexamples to show that the Wasserstein barycenter of elliptically contoured distributions does not have to be elliptically contoured.
\end{abstract}

% keywords can be removed

\section{Introduction}
Optimal transport has attracted considerable attention in recent years. It has been a powerful tool in mathematics \cite{EvansL.C1999DEMf,Villani2009OTOa}, physics \cite{HakerSteven2004Omtf}, engineering \cite{6502714}, and data science \cite{korotin2023neural}. This theory established a metric structure on the space of probability distributions. It was originally considered by Monge \cite{monge1781memoire}, and studied deeply after the relaxation by Kantorovich \cite{KantorovichL.V.2006Otto}. When the two probability measures are absolutely continuous with respect to the Lebesgue measure, the formulations given by Monge and Kantorovich are equivalent \cite{BrenierYann1991Pfam}. The optimal transport theory also provides a natural interpolation between probability distributions \cite{McCannRobertJ.1997ACPf}. The interpolation can be generalized to more than two distributions, which is known as Wasserstein barycenter \cite{AguehMartial2011Bitw}. The optimal transport map and Wasserstein barycenter play an important role in practice \cite{SolomonJustin2015CWdE}. However, only in a few cases, we can get the specific property of the optimal transport map and Wasserstein barycenter. For instance, Gaussian distributions exhibit closed-form solutions: their optimal transport maps are affine transformations, and Wasserstein barycenters preserve Gaussianity \cite{PeyreGabriel2019Cot}. The result can be easily generalized to elliptically contoured distributions from the same family \cite{GelbrichMatthias1990OaFf}.

Gaussian mixture models (GMMs) have been popular in applications. They are widely used to represent distributions of real datasets \cite{DeledalleCharlesAlban2018Idwg,GalerneBruno2017SOTi,TeodoroAfonsoM.A.M.2015Sida}. Recently, a relaxed formula of optimal transport for GMMs has been proposed \cite{DelonJulie2020Awdi,ChenYukun2020AWDa,ChenYongxin2019OTfG}. The theoretical frame is based on the fact that the McCann interpolation and Wasserstein barycenter of Gaussian distributions are still Gaussian. The generic mixtures whose components are elliptically contoured distributions have been investigated in \cite{dusson2023wassersteintypemetricgenericmixture}. However, in their studies, the elliptically contoured distributions always share a same generator function. Theoretically, while approximating two or more unknown distributions, using mixtures in different families may yield superior accuracy than using the same type mixtures. To generalize their results to more than two families of mixtures, it is essential to study the property of the optimal transport and Wasserstein barycenter for elliptically contoured distributions with different generator functions. However, to our knowledge, the corresponding results for elliptically contoured distributions with different generator functions have not been discussed in previous study.

In this paper, we consider a specialized class of elliptically contoured distributions by degenerating them to radially contoured distributions. Such distributions were also studied in \cite{CuestaalbertosJ.A.1993OCoM}. We derive the formulas of the Monge map and the Wasserstein distance between two radially contoured distributions with different generator functions and prove that the Wasserstein barycenter of radially contoured distributions remains radially contoured. However, similar results do not hold true for general elliptically contoured distributions. We give two numerical counterexamples to show that the Wasserstein barycenter of elliptically contoured distributions does not have to be elliptically contoured.

This paper is organized as follows. Section \ref{sec2} presents related results in optimal transport theory and Wasserstein barycenter. In section \ref{sec3}, we study the optimal transport for radially contoured distributions. Section \ref{sec4} provides a rigorous proof of the structural preservation property for Wasserstein barycenters of radially contoured distributions. In section \ref{sec5}, two numerical counterexamples are given to show that the Wasserstein barycenter of elliptically contoured distributions does not have to be elliptically contoured. Finally, in section \ref{sec6}, we conclude this paper and discuss about the further study.

\textbf{Notation} We define some of the notation that will be used in this paper.
\begin{itemize}
	\item $\R_+ = [0,+\infty)$ denotes the set of nonnegative real numbers.
	\item $\norm{x}$ denotes the Euclidean norm of $x \in \R^d$ and $\inner{x}{y}$ denotes the Euclidean inner product of $x$ and $y$ in $\R^d$.
	\item $\mathbb{S}^{d-1}$ denotes the $d-1$ sphere in $\R^d$, i.e. $\mathbb{S}^{d-1} = \{x \in \R^d: \norm{x} = 1\}$.
	\item $\Gamma_N = \left\{ \pi = (\pi_j)_{j=1}^N \in \R_{+}^N : \sum_{j=1}^N \pi_j = 1 \right\}$ denotes the probability simplex.
	\item $I_d$ denotes the $d \times d$ identity matrix.
	%	\item When $a_i$ is a finite sequence of $K$ elements, we denotes its elements as $a_i^1,\cdots,a_i^K$. 
	\item  $P_0$ and $P_1$ denote the projections on $\R^d \times \R^d$ such that $P_0(x,y) = x, P_1(x,y) = y$.
	\item If $\mu$ is a positive measure on a space $\mathcal{X}$ and $T: \mathcal{X} \to \mathcal{Y}$ is a measurable map, $T_\# \mu$ stands for the push-forward measure of $\mu$ by $T$, i.e. the measure on $\mathcal{Y}$ such that for all measurable set $A \subset \mathcal{Y}$, $(T_\# \mu)(A) = \mu (T^{-1}(A))$.
	%	\item If $\mu$ is a probability measure on a space $\mathcal{X}$, $\bar{\mu}$ denotes centralization of $\mu$.
\end{itemize}

\section{Review of related results}\label{sec2}

Let $d \geq 1$ be an integer and $\mathcal{P}(\R^d)$ denote the set of probability measures on $\R^d$. For $p \geq 1$, the Wasserstein space $\mathcal{P}_p(\R^d)$ is defined as the set of probability measures with finite moment of order $p$, such that
\begin{equation*}
	\mathcal{P}_p(\R^d) = \left\{ \mu \in \mathcal{P}(\R^d): \int_{\R^d} \norm{x}^p \mathrm{d} \mu(x) < + \infty \right\}.
\end{equation*}
From now on, we will focus on the case $p=2$.

Let $\mu_0, \mu_1 \in \mathcal{P}_2(\R^d)$, the Monge formulation of optimal transport is
\begin{equation}\label{Monge}
	\inf_{T: T_\# \mu_0 = \mu_1} \int_{\R^d} \norm{x - T(x)}^2 \mathrm{d} \mu_0(x),
\end{equation}
where $T$ takes over all measurable maps from $\R^d$ to $\R^d$. The optimal $T$ is called Monge map transporting $\mu_0$ to $\mu_1$. Define $\Pi(\mu_0,\mu_1) \subset \mathcal{P}_2(\R^d \times \R^d)$ as the subset of probability distributions $\gamma$ on $\R^d \times \R^d$ with marginals $\mu_0$ and $\mu_1$, such that $(P_0)_\# \gamma = \mu_0$ and $(P_1)_\# \gamma = \mu_1$. 
%More over, if $\mu_0 = \pi_0 \in \Gamma_{K_0}$ and $\mu_1 = \pi_1 \in \Gamma_{K_1}$, $\Pi(\pi_0,\pi_1)$ is represented as
%\begin{equation*}
%	\Pi(\pi_0,\pi_1) = \left\{w \in \R_+^{K_0 \times K_1}: \sum_{l=1}^{K_2} w_{kl} = \pi_0^k, \, \sum_{k=1}^{K_1} w_{kl} = \pi_1^l\right\}.
%\end{equation*}
Kantorovich gives a relaxed form of (\ref{Monge}) by
\begin{equation}\label{W2}
	W_2(\mu_0, \mu_1)^2 = \inf_{\gamma \in \Pi(\mu_0, \mu_1)} \int_{\R^d\times\R^d} \norm{x-y}^2 \mathrm{d} \gamma(x,y).
\end{equation}
$W_2(\cdot,\cdot)$ actually defines  a distance on $\mathcal{P}_2(\R^d)$, which is known as the Wasserstein distance. In general case, (\ref{Monge}) and (\ref{W2}) are not equivalent. However, thanks to Brenier, Theorem \ref{thm:Brenier} gives the equivalence of the two problems when $\mu_0$ has a density, which is illustrated in the following.
\begin{theorem}[Brenier \cite{BrenierYann1991Pfam}, \cite{PeyreGabriel2019Cot}]\label{thm:Brenier}
	Given two measures $\mu_0$ and $\mu_1$ on $\R^d$. Suppose $\mu_0$ has a density $\rho$ with respect to the Lebesgue measure, then the optimal $\gamma$ in $(\ref{W2})$ is unique and supported on the graph $(x,T(x))$ of a Monge map $T: \R^d \to \R^d$. It means $\gamma = (\operatorname{Id}, T)_\# \mu_0$. Furthermore, $T$ is uniquely defined as the gradient of a convex function $\phi$, $T(x) = \nabla \phi(x)$, where $\phi$ is the unique (up to an additive constant) convex function such that $(\nabla \phi)_\# \mu_0 = \mu_1$.
\end{theorem}

A nice property of $W_2$ is one can factor out translations. Let $m_i = \int_{\R^d} x \mathrm{d} \mu_i (x)$ be the mean of $\mu_i$, and $T_\tau: x \mapsto x - \tau$ be that translation operator, one has
\begin{equation*}
	W_2(\mu_0,\mu_1)^2 = W_2(\bar{\mu}_0,\bar{\mu}_1)^2 + \norm{m_0 - m_1}^2, \quad \bar{\mu}_i = (T_{m_i})_\#\mu_i, \quad i = 0,1.
\end{equation*}
If $T$ is a Monge map transporting $\mu_0$ to $\mu_1$, the path $(\mu_t)_{t\in[0,1]}$ given by
\begin{equation}\label{geodesic curve}
	\mu_t = (T_t)_\# \mu_0, \quad T_t(x) = (1-t)x + tT(x), \quad \forall \, t \in [0,1],
\end{equation}
defines a constant speed geodesic curve on $(\mathcal{P}_2(\R^d),W_2)$, i.e.
\begin{equation}\label{geodesic property}
	W_2(\mu_s, \mu_t) = (t-s) W_2(\mu_0,\mu_1), \quad 0 \leq s \leq t \leq 1.
\end{equation}
This path is called McCann interpolation \cite{McCannRobertJ.1997ACPf} between $\mu_0$ and $\mu_1$. It satisfies
\begin{equation}\label{mu_t}
	\mu_t \in \mathrm{argmin}_{\mu \in \mathcal{P}_2(\R^d)} (1-t)W_2(\mu_0, \mu)^2 + tW_2(\mu_1,\mu)^2.
\end{equation}

The McCann interpolation can be generalized to more than two distributions, which is known as Wasserstein barycenter. Given an integer $N \geq 2$, $N$ probability measures $\mu_1, \mu_2, \cdots, \mu_N$ and weights $\lambda = (\lambda_1, \lambda_2, \cdots, \lambda_N) \in \Gamma_N$, the associated Wasserstein barycenter problem is:
\begin{equation}\label{Barycenter Problem}
	\inf_{\mu \in \mathcal{P}_2(\R^d)} \sum_{j=1}^N \lambda_j W_2(\mu_j,\mu)^2.
\end{equation}
A solution of the previous problem is called the Wasserstein barycenter of probability measures $\{\mu_j\}_{j=1}^N$ with weights $\{\lambda_j\}_{j=1}^N$. The existence and uniqueness have been deeply studied in \cite{AguehMartial2011Bitw}. Here we quote Theorem \ref{thm:barycenter characterization} below which gives a characterization of the Wasserstein barycenter.

\begin{theorem}[Characterization of barycenters \cite{AguehMartial2011Bitw}]\label{thm:barycenter characterization}
	Assume that $\{\mu_j\}_{j=1}^N$ are distributions on $\R^d$, and each one vanishes on small sets ($\mu_j(A) = 0$ for all Borel sets $A$ of $\R^d$, having Haussdorff dimension less than $d-1$). Let $\mu \in \mathcal{P}_2(\R^d)$. Then the following conditions are equivalent:
	
	1. $\mu$ uniquely solves $(\ref{Barycenter Problem})$.
	
	2. $\mu = (\nabla \phi_j)_\# \mu_j$ for every $j$, where $\phi_j$ is defined by
	\begin{equation}\label{phi}
		\lambda_j \phi_j := \frac{\lambda_j}{2} \norm{x}^2 - S_{\lambda_j} f_j(x).
	\end{equation}
	Here $S_\lambda f(x)$ is defined by
	\begin{equation}
		S_\lambda f(x) := \inf_{y \in \R^d} \left\{ \frac{\lambda}{2} \norm{x-y}^2 - f(y) \right\}, \quad \forall \, x \in \R^d,
	\end{equation}
	and $\{f_j\}_{j=1}^N$ are the solutions of the following problem:
	\begin{equation}\label{dual}
		\sup \left\{ F(f_1, f_2, \cdots, f_N) = \sum_{j=1}^N \int_{\R^d} S_{\lambda_j} f_j(x) \mathrm{d} \mu_j(x): \sum_{j=1}^{N} f_j = 0, \, f_j \in X, \forall \, j  \right\},
	\end{equation}
	\begin{equation*}
		X := \left\{ f(x) \in C(\R^d): \frac{f(x)}{1 + \norm{x}^2} \, \text{is bounded} \right\}.
	\end{equation*}
	
	3. There exist convex potentials $\{\psi_j\}_{j=1}^N$ such that $\nabla \psi_j$ is the Monge map transporting $\mu_j$ to $\mu$, and a constant $C$ such that
	\begin{equation}
		\sum_{j=1}^{N} \lambda_j \psi_j^* (y) \leq \frac{\norm{y}^2}{2} + C, \quad \forall \, y \in \R^d, \, \mu \text{-a.e.} 
	\end{equation}
	Here $\psi_j^*$ denotes the Legendre transformation of $\psi_j$:
	\begin{equation}
		\psi_j^*(y) = \sup_{x \in \R^d} \left\{ \inner{x}{y} - \psi_j(x) \right\}.
	\end{equation}
\end{theorem}
\begin{remark}
	A probability measure with a continuous density vanishes on small sets.
\end{remark}

\section{Optimal transport for radially contoured distributions}
\label{sec3}

\begin{definition}[Radially contoured distribution]\label{def Radial}
	A probability measure $\mu$ on $\R^d$ is a radially contoured distribution centered on $m$, if $\mu$ is a probability measure and
	\begin{equation*}
		\mathrm{d}\mu(x) = \frac{1}{Z}\rho\left(\frac{\norm{x-m}}{c}\right) \dx,
	\end{equation*}
	where $\rho$ is a nonnegative measurable function on $\R_+$ and $Z$ is the normalizer that ensures $\int_{\R^d} \mathrm{d} \mu(x) = 1$. We call $m$ the center and $\rho$ the generator function of $\mu$. We denote $\mu = R_d(m,c,\rho)$. If $Z=c=1$, we simply denote $\mu = R_d(m,\rho)$. If $m = 0$, we will call $\mu$ a radially distribution. Moreover, we assume $\mu \in \mathcal{P}_2(\R^d)$ in this paper.
\end{definition}

\begin{remark}
	Notice that for all $\mu = R_d(m,c,\rho)$, we can define $\tilde{\rho}(\cdot) = Z^{-1}\rho(\cdot/c)$. This means that we can always assume $Z=c=1$. i.e. $\mu = R_d(m,\tilde{\rho})$. The notation $R_d(m,c,\rho)$ will be used in specific examples, and $R_d(m,\rho)$ will be used in theorems about optimal transport and Wasserstein barycenter of radially contoured distributions with different generator functions.
\end{remark}

\begin{remark}
	This definition can be viewed as the degeneration case of elliptically contoured distribution when the positive definite matrix is a scalar matrix. In most time, we assume the generator function $\rho$ is continuous.
\end{remark}

By straightforward calculation, we can get the normalizer, expectation and covariance of a radially contoured distribution as follows.

\begin{proposition}\label{prop1}
	Let $X \sim \mu = R_d(m, c, \rho) \in \mathcal{P}_2(\R^d)$, then
	\begin{equation*}
		Z = c^d |\mathbb{S}^{d-1}| \int_0^{+\infty} r^{d-1} \rho(r) \dr, \quad  E[X] = m, \quad \operatorname{Cov}[X] = \frac{c^2}{d} \frac{\int_0^{+\infty}r^{d+1}\rho(r)\dr}{\int_0^{+\infty}r^{d-1}\rho(r)\dr}I_d,
	\end{equation*}
	where $|\mathbb{S}^{d-1}| = 2 \pi^{d/2} / \Gamma(d/2)$ denotes the volume of $\mathbb{S}^{d-1}$.
\end{proposition}

The main purpose of this section is to study the optimal transport for radially contoured distributions here. We give formulas of the Monge map, the McCann interpolation and Wasserstein distance for radially contoured distributions in Theorem \ref{thm:r_ot}. To prove the results, we first give two lemmas.

\begin{lemma}\label{lem C}
	Let $\mu_0, \mu_1$ are probability measures on $\R_+$  with continuous density functions $\rho_0$ and $\rho_1$. Then there exists a nondecreasing function $C: \R_+ \to \R_+ \cup \{+\infty\}$ that transports $\mu_0$ to $\mu_1$, and
	\begin{equation}\label{c mu0 mu1}
		\mu_0([0,R]) = \mu_1([0,C(R)]), \quad \forall \, R \in \R_+.
	\end{equation}
\end{lemma}

\textbf{Proof.} This lemma can be proved in a similar way to prove Theorem 2.5 in \cite{SantambrogioFilippo2015OTfA}.

\begin{remark}
	The existence of the nondecreasing transport map $C$ only relies on the atomless property of $\mu_0$. However, (\ref{c mu0 mu1}) also requires the atomless property of $\mu_1$. Rather than the construction with the pseudoinverse of the cumulative distribution function, (\ref{c mu0 mu1}) gives a more intuitive description of $C$.
\end{remark}

\begin{lemma}\label{lem convex}
	Let $f: \R_+ \to \R_+ \cup \{ + \infty\}$ be a nondecreasing nonnegative function, then for all $m \in \R^d$,
	\begin{equation*}
		u(x) = \int_0^{\norm{x - m}} f(r) \dr, \quad \forall \, x \in \R^d,
	\end{equation*}
	is a convex function on $\R^d$.
\end{lemma}

\textbf{Proof.} We only prove the case $m=0$. Let $F(R) = \int_0^R f(r) \dr$. Since $f$ is nondecreasing, $F$ is convex on $\R_+$. For all $x_1, x_2 \in \R^d$ and $\lambda \in [0,1]$, we have
\begin{align*}
		u(\lambda x_1 + (1-\lambda) x_2) & = F(\norm{\lambda x_1 + (1-\lambda) x_2}) \\ & \leq F(\lambda \norm{x_1} + (1-\lambda) \norm{x_2}) \\ & \leq \lambda F(\norm{x_1}) + (1 - \lambda) F(\norm{x_2}) \\ & = \lambda u(x_1) + (1 - \lambda) u(x_2).
\end{align*}

\begin{theorem}\label{thm:r_ot}
	Let $\mu_0 = R_d(m_0,\rho_0), \, \mu_1 = R_d(m_1,\rho_1)$ be two radially contoured distributions on $\R^d$ with continuous generator functions $\rho_0$ and $\rho_1$. Then there exists a function $C : \R_+ \to \R_+ \cup \{+\infty\}$ satisfying
	\begin{equation}\label{C}
		\int_0^R r^{d-1} \rho_0(r) \dr = \int_0^{C(R)} r^{d-1} \rho_1(r) \dr,
	\end{equation}
	and the Monge map transporting $\mu_0$ to $\mu_1$ is given by
	\begin{equation*}
		T(x) = \frac{C\left(\norm{x-m_0}\right)}{\norm{x-m_0}}(x-m_0) + m_1.
	\end{equation*}
	The McCann interpolation is given by $\mu_t = R_d(m_t,\rho_t)$, where $m_t = (1-t)m_0 + t m_1$ and $\rho_t$ satisfies
	\begin{equation}\label{rho_t}
		\left(C_t\right)_\# R_1\left(0,r^{d-1}\rho_0(r)\right) = R_1\left(0,r^{d-1}\rho_t(r)\right),
	\end{equation}
	%	\begin{equation}\label{rho_t}
		%		\rho_t\left(C_t(r)\right)C_t(r)^{d-1} \mathrm{d} C_t(r) = \rho_0(r) r^{d-1} \dr.
		%	\end{equation}
	where $C_t(r) = (1-t)r + tC(r)$. The Wasserstein distance is given by
	\begin{align*}
		W_2(\mu_0,\mu_1)^2 & = \norm{m_0 - m_1}^2 + W_2(R_1(0,r^{d-1}\rho_0(r)),R_1(0,r^{d-1}\rho_1(r))) \\ & = \norm{m_0 - m_1}^2 + |\mathbb{S}^{d-1}| \int_0^{+\infty} \left(C(r) - r\right)^2 \rho_0(r) r^{d-1} \dr.
	\end{align*}
\end{theorem}

\textbf{Proof.} We have
	\begin{equation*}
		\int_{\R^d} \rho_0(\norm{x}) \dx = \int_{\R^d} \rho_1(\norm{x}) \dx = 1.
	\end{equation*}
	After spherical coordinate transformation, it turns out that
	\begin{equation*}
		\int_0^{+\infty} r^{d-1} \rho_0(r) \dr = \int_0^{+\infty} r^{d-1} \rho_1(r) \dr = \frac{1}{|\mathbb{S}^{d-1}|}.
	\end{equation*}
	Let $\tilde{\mu}_0, \tilde{\mu}_1$ be two probability measures on $\R$ such that
	\begin{equation*}
		\mathrm{d} \tilde{\mu}_i (r) = |\mathbb{S}^{d-1}| r^{d-1} \rho_i(r) \mathrm{d} r, \quad i = 0, 1 .
	\end{equation*}
	By Lemma \ref{lem C}, there exists a nondecreasing function $C: \R_+ \to \R_+\cup \{+\infty\}$ such that $C_\# \tilde{\mu}_0 = \tilde{\mu}_1$. It is equivalent to (\ref{C}).
	
	%	Because
	%	\begin{equation*}
		%		C(R) = G^{-1}(F(R)),
		%	\end{equation*}
	%	we get that $C$ is differentiable.
	Consider $u(x) = \int_0^{\norm{x-m_0}} C(r) \dr + \inner{m_1}{x}$. By Lemma \ref{lem convex}, $u$ is a convex function. Moreover,
	\begin{equation*}
		\nabla u(x) = C(\norm{x-m_0}) \frac{x-m_0}{\norm{x-m_0}} + m_1 = T(x),
	\end{equation*}
	%	Taking derivative on both side of (\ref{C}), we get
	%	\begin{equation*}
		%		r^{d-1} \rho_0(r) = C(r)^{d-1}\rho_1(C(r))C^\prime(r).
		%	\end{equation*}
	and
	\begin{equation*}\label{analogous}
		\begin{aligned}
			\int_{\R^d} h(x) \mathrm{d} T_\# \mu_0(x) & = \int_{\R^d} h\left(T(x)\right) \rho_0\left(\norm{x-m_0}\right) \dx \\
			& = \int_{\mathbb{S}^{d-1}} \dSz \int_0^{+\infty} h\left(T(r \cdot z + m_0)\right) r^{d-1} \rho_0(r) \dr \\
			& = \int_{\mathbb{S}^{d-1}} \frac{1}{|\mathbb{S}^{d-1}|} \dSz \int_0^{+\infty} h(C(r) \cdot z + m_1) \mathrm{d} \tilde{\mu}_0 (r) \\
			& = \int_{\mathbb{S}^{d-1}} \frac{1}{|\mathbb{S}^{d-1}|} \dSz \int_0^{+\infty} h(r \cdot z + m_1) \mathrm{d} \tilde{\mu}_1 (r) \\
			& = \int_{\mathbb{S}^{d-1}} \dSz \int_0^{+\infty} h(r \cdot z + m_1) r^{d-1} \rho_1(r) \dr \\
			& = \int_{\R^d} h(x) \mathrm{d} \mu_1(x), \quad \forall \, h \in C({\R^d}).
		\end{aligned}
	\end{equation*}
	The fourth equation is given by $C_\# \tilde{\mu}_0 = \tilde{\mu}_1$ which means $T_\# \mu_0 = \mu_1$. By Theorem \ref{thm:Brenier}, $T$ is the Monge map. 
	
	Taking $T_t(x) = (1-t)x + t T(x)$, $C_t(r) = (1-t)r + t C(r)$ and $\rho_t$ in (\ref{rho_t}), since $\mu_t = (T_t)_\# \mu_0$, we have
	\begin{align*}
		\int_{\R^d} h(x) \mathrm{d} (T_t)_\# \mu_0(x) & = \int_{\mathbb{S}^{d-1}} \frac{1}{|\mathbb{S}^{d-1}|} \dSz \int_0^{+\infty} h\left(T_t(r \cdot z + m_0)\right)r^{d-1} \rho_0(r) \dr \\
		& = \int_{\mathbb{S}^{d-1}} \frac{1}{|\mathbb{S}^{d-1}|} \dSz \int_0^{+\infty} h\left(C_t(r) \cdot z + m_t\right)r^{d-1} \rho_0(r) \dr \\
		& = \int_{\mathbb{S}^{d-1}} \frac{1}{|\mathbb{S}^{d-1}|} \dSz \int_0^{+\infty} h(r \cdot z + m_t)r^{d-1} \rho_t(r) \dr.
	\end{align*}
	This means that $\mu_t = R_d(m_t, \rho_t)$.
	
	Let $\bar{\mu}_i = (T_{m_i})_\#\mu_i$, then
	\begin{equation*}
		W_2(\mu_0,\mu_1)^2 = \norm{m_0 - m_1}^2 + W_2(\bar{\mu}_0, \bar{\mu}_1)^2.
	\end{equation*}
	
	By the previous proof, $C$ is the Monge map that transports $R_1(0,r^{d-1}\rho_0(r))$ to $R_1(0,r^{d-1}\rho_1(r))$, and the Monge map transporting $\bar{\mu}_0$ to $\bar{\mu}_1$ is
	\begin{equation*}
		\hat{T}(x) = C(\norm{x}) \frac{x}{\norm{x}}. 
	\end{equation*}
	Therefore, we have
	\begin{align*}
		%		W_2(\bar{\mu}_0, \bar{\mu}_1)^2 = \int_{\R^d} \norm{x - \hat  {T}(x)}^2 \rho_0(\norm{x}) \dx = |\mathbb{S}^{d-1}| \int_0^{+\infty} (r - C(r))^2 \rho_0(r) r^{d-1} \dr.
		W_2(\bar{\mu}_0, \bar{\mu}_1)^2 & = \int_{\R^d} \norm{x-T(x)}^2 \mathrm{d} \bar{\mu}_0(x) = \int_{\mathbb{S}^{d-1}} \dSz \int_0^{+\infty} \left(r-C(r)\right)^2 r^{d-1} \rho_0(r) \dr \\
		& = |\mathbb{S}^{d-1}| \int_0^{+\infty} \left(C(r) - r\right)^2 r^{d-1} \rho_0(r)\dr \\
		& = |S^{d-1}| W_2\left(R_1\left(0,r^{d-1}\rho_0(r)\right), R_1\left(0,r^{d-1} \rho_1(r)\right)\right)^2.
	\end{align*}

When $\mu_0$ and $\mu_1$ share the same generator function, the Monge map and Wasserstein distance have simple formulas.

\begin{example}\label{eg1}
	Let $\mu_0 = R_d(m_0, c_0, \rho)$, $\mu_1 = R_d(m_1, c_1, \rho)$, then the Monge map from $\mu_0$ to $\mu_1$ is
	\begin{equation*}
		T(x) = \frac{c_1}{c_0}(x-m_0) + m_1.
	\end{equation*}
	The McCann interpolation is $\mu_t = R_d(m_t,c_t,g)$, where
	\begin{equation*}
		m_t = (1-t)m_0 + t m_1, \quad c_t = (1-t) c_0 + t c_1, \quad 0 \leq t \leq 1,
	\end{equation*}
	and
	\begin{equation}
		W_2(\mu_0,\mu_1)^2 = \norm{m_0 - m_1}^2 + (c_0 - c_1)^2 \frac{\int_0^{+\infty} r^{d+1} \rho(r)\dr}{\int_0^{+\infty} r^{d-1} \rho(r) \dr}.
	\end{equation}
\end{example}
\textbf{Proof.} By viewing the radially counted distributions as degenerated elliptically contoured distributions, the proof of this theorem can be completed by synthesizing the conclusion of \cite{GomezEusebio2003Asoc,GelbrichMatthias1990OaFf,PeyreGabriel2019Cot}.
	
	Here, we provide the following proof using Theorem \ref{thm:r_ot}. Let $\rho_0(r) = Z_0^{-1}\rho(r/c_0)$, $\rho_1(r) = Z_1^{-1}\rho(r/c_1)$, where $Z_0$ and $Z_1$ are normalizers. It is easy to verify that
	\begin{equation*}
		\rho_1(r) = \left(\frac{c_0}{c_1}\right)^d \rho_0\left(\frac{c_0}{c_1} r\right).
	\end{equation*}
	Then we find $C(r) = c_1 r / c$ is the solution of (\ref{C}). Therefore, the Monge map is
	\begin{equation*}
		T(x) = \frac{c_1}{c_0}(x-m_0) + m_1.
	\end{equation*}
	Since $C_t(r) = (1-t)r + (tc_1/c_0)r$, we have $C_t^{-1}(r) = c_0 r / c_t$. Therefore
	\begin{equation*}
		\rho_t(r) = \left(\frac{c_0}{c_t}\right)^d \rho_0\left(\frac{c_0}{c_t} r\right) = Z_0^{-1} \left(\frac{c_0}{c_t}\right)^d \rho\left(\frac{r}{c_t} \right),
	\end{equation*}
	which means $\mu_t = R_d(m_t,c_t,\rho)$. At last, we calculate that
	\begin{align*}
		& W_2\left(R_1\left(0,r^{d-1}\rho_0(r)\right), R_1\left(0,r^{d-1} \rho_1(r)\right)\right)^2 = \frac{1}{c_0^2} (c_0 - c_1)^2 \int_{\R^d} \norm{x}^2 \rho_0\left(\norm{x}\right) \dx \\
		= & \frac{1}{c_0^2} (c_0 - c_1)^2 \operatorname{trace}\left(\operatorname{Cov}\left(R_d(0,\rho_0)\right)\right) = \frac{1}{c_0^2} (c_0 - c_1)^2 \operatorname{trace}\left(\operatorname{Cov}\left(R_d(0,c_0,\rho)\right)\right) \\
		= & (c_0 - c_1)^2 \frac{\int_0^{+\infty}r^{d+1}\rho(r)\dr}{\int_0^{+\infty}r^{d-1}\rho(r)\dr}.
	\end{align*}

\section{Wasserstein barycenter of radially contoured distributions}\label{sec4}

In this section, we will study the Wasserstein barycenter of radially contoured distributions. In Theorem \ref{thm:r_ot}, we get that the McCann interpolation of two radially contoured distributions is still radially contoured. Here, we show that while generalizing to more than two radially contoured distributions, similar result will hold. Their barycenter is a generalized radially contoured distribution.
%  It will be illustrated in Theorem \ref{thm:radial barycenter}.

\begin{definition}(Generalized radially contoured distribution)\label{def Gradial}
	A probability measure $\mu$ on $\R^d$ is a generalized radially contoured distribution centered on $m$ if there exists a measure $\nu$ on $\R_+$, such that for all $h(x) \in C(\R^d)$,
	\begin{equation}\label{irregular}
	\int_{\R^d} h(x) \mathrm{d} \mu(x) = \int_{\mathbb{S}^{d-1}} \int_0^{+\infty} h(r \cdot z + m) \mathrm{d} \nu(r) \dSz.
	\end{equation}
	If $m = 0$, we will call $\mu$ a generalized radial distribution.
\end{definition}

\begin{remark}
	Compared to Definition \ref{def Radial}, a singular part is added in Definition \ref{def Gradial}. To see this, let $F(r) = \nu([0,r])$ be the cumulative distribution function. We can write $F = F^a + F^j$, where $F^a$ is an absolutely continuous function and $F^j$ is a jumping function with countable right continuous jumping points. The first part is associated with a radially contoured distribution, and the second part is the cumulative distribution function some Dirac measures. The generalized derivative of $F$ can be viewed as the generator function of the generalized radially contoured distribution.
\end{remark}

\begin{remark}
	A generalized radially contoured distribution vanishes on small sets.
\end{remark}

The main tools to prove our theorem are Theorem \ref{thm:Brenier} and Theorem \ref{thm:barycenter characterization}. Moreover, we need the following $3$ lemmas. 
\begin{lemma}\label{lem1}
	Let $\phi$ be a radial convex function, and $\mu = R_d(0,\rho)$ be a radial distribution. Then $(\nabla \phi)_\# \mu$ is a generalized radial distribution.
\end{lemma}

\textbf{Proof.} Since $\phi$ is radial and convex, $\nabla \phi(x) = T(\norm{x}) \cdot n_x$ almost everywhere, and $T$ is a nondecreasing function on $[0, +\infty)$. We define
\begin{equation*}
	\mathrm{d} \mu_{rad}(r) := r^{d-1} \rho(r) \dr, \quad r \geq 0,
\end{equation*}
and $\nu = T_\# \mu_{rad}$. For all $h \in C(\R^d)$, we have
\begin{align*}
	& \int_{\R^d} h(x) \mathrm{d} (\nabla \phi)_\# \mu(x) = \int_{\R^d} h(\nabla \phi(x)) \rho(\norm{x}) \dx \\
	= & \int_{\mathbb{S}^{d-1}} \dSz \int_0^{+\infty} h(T(r) \cdot z) r^{d-1} \rho(r) \dr = \int_{\mathbb{S}^{d-1}} \dSz \int_0^{+\infty} h(T(r) \cdot z) \mathrm{d} \mu_{rad} (r) \\
	= & \int_{\mathbb{S}^{d-1}} \dSz \int_0^{+\infty} h(r \cdot z) \mathrm{d} \nu (r).
\end{align*}

\begin{lemma}\label{lem2}
	Let $f(\cdot) = g(\norm{\cdot})$ be a radial function on $\R^d$. Then $S_\lambda f$ is still a radial function on $\R^d$.
\end{lemma}

\textbf{Proof.} For all $x \in \R^d$, we write $x = R_x n_x$, where $R_x = \norm{x}$. For all positive $R$ and $n \in \mathbb{S}^{d-1}$,
	\begin{equation*}
		\norm{R_x \cdot n_x - R \cdot n_x} = |R_x - R| \leq \norm{R_x \cdot n_x - R \cdot n}. 
	\end{equation*}
	Therefore
	\begin{equation*}
		\frac{\lambda}{2} \norm{R_x \cdot n_x - R \cdot n_x} - f(R \cdot n_x) \leq \frac{\lambda}{2} \norm{R_x \cdot n_x - R \cdot n} - f(R \cdot n).
	\end{equation*}
	So we can get that
	\begin{align*}
		S_\lambda f(x) & = \inf_{y \in \R^d} \left\{ \frac{\lambda}{2}\norm{x-y}^2 - f(y) \right\} \\
		& = \inf_{R \geq 0, n \in \mathbb{S}^{d-1}} \left\{ \frac{\lambda}{2}\norm{x-R \cdot n}^2 - f(R \cdot n) \right\} \\
		& = \inf_{R \geq 0} \left\{ \frac{\lambda}{2}\norm{x-R \cdot n_x}^2 - f(R \cdot n_x) \right\} \\
		& = \inf_{R \geq 0} \left\{ \frac{\lambda}{2}(\norm{x}-R)^2 - g(R) \right\}.
	\end{align*}
	For $x_1, x_2 \in \R^d$ satisfied with $\norm{x_1} = \norm{x_2}$, we conclude that
	\begin{equation*}
		S_{\lambda} f(x_1) = \inf_{R \geq 0} \left\{ \frac{\lambda}{2}(\norm{x_1}-R)^2 - g(R) \right\} = \inf_{R \geq 0} \left\{ \frac{\lambda}{2}(\norm{x_2}-R)^2 - g(R) \right\} = S_{\lambda} f(x_2).
	\end{equation*}

\begin{lemma}\label{lem3}
	Given a function $f$ and a radial distribution $\mu = R_d(0,\rho)$ on $\R^d$, we define
	\begin{equation}\label{radial f}
		\bar{f}(x) = g_f(\norm{x}) = \begin{cases}
			\frac{1}{|\partial B_{\norm{x}}|}\int_{\partial B_{\norm{x}}} f(x) \dx, & \norm{x} > 0, \\
			f(0),  & x = 0,
		\end{cases}
	\end{equation}
	as the radialization of $f$. Then
	\begin{equation}
		\int_{\R^d} S_\lambda f(x) \mathrm{d} \mu(x) \leq \int_{\R^d} S_\lambda \bar{f}(x) \mathrm{d} \mu(x).
	\end{equation}
\end{lemma}

\textbf{Proof.} We first simplify the right-hand side of the inequality.
	\begin{align*}
		& \int_{\R^d} S_\lambda \bar{f} (x) \mathrm{d} \mu (x) \\ = & \int_0^{+\infty}r^{d-1}\rho(r) \dr \int_{\mathbb{S}^{d-1}} \inf_{y \in \R^d} \left\{ \frac{\lambda}{2} \norm{r \cdot z - y}^2 -f(y) \right\} \dSz \\
		= & \int_0^{+\infty}r^{d-1}\rho(r) \dr \int_{\mathbb{S}^{d-1}} \inf_{R \geq 0, n \in \mathbb{S}^{d-1}} \left\{ \frac{\lambda}{2} \norm{r \cdot z - R \cdot n} - g_f(R) \right\} \dSz \\
		= & \int_0^{+\infty}r^{d-1}\rho(r) \dr \int_{\mathbb{S}^{d-1}} \inf_{R \geq 0} \left\{ \frac{\lambda}{2} (r - R)^2 - g_f(R) \right\} \dSz \label{eq:a}\\
		= & \int_0^{+\infty} |\mathbb{S}^{d-1}| \inf_{R \geq 0}\left\{ \frac{\lambda}{2} (r-R)^2 - g_f(R)\right\} r^{d-1}\rho(r) \dr.
	\end{align*}
	The third equation is an application of the proof of Lemma \ref{lem2}. Therefore,
	\begin{align*}
		& \int_{\R^d} S_\lambda f (x) \mathrm{d} \mu(x) \\ = & \int_0^{+\infty}r^{d-1}\rho(r) \dr \int_{\mathbb{S}^{d-1}} \inf_{y \in \R^d} \left\{ \frac{\lambda}{2} \norm{r \cdot z - y}^2 -f(y) \right\} \dSz \\
		= & \int_0^{+\infty}r^{d-1}\rho(r) \dr \int_{\mathbb{S}^{d-1}} \inf_{\substack{R \geq 0, n \in \mathbb{S}^{d-1}}} \left\{ \frac{\lambda}{2} \norm{r \cdot z - R \cdot n}^2 - f(R \cdot n)\right\} \dSz \\
		\leq & \int_0^{+\infty}r^{d-1}\rho(r) \dr \int_{\mathbb{S}^{d-1}} \inf_{R \geq 0} \left\{ \frac{\lambda}{2}\norm{r \cdot z - R \cdot z}^2 - f(R \cdot z) \right\} \dSz \\
		\leq & \int_0^{+\infty}r^{d-1}\rho(r) \dr \inf_{R \geq 0} \int_{\mathbb{S}^{d-1}} \left\{ \frac{\lambda}{2} (r-R)^2 - f(R \cdot z) \right\} \dSz \\
		= & \int_0^{+\infty} |\mathbb{S}^{d-1}| \inf_{R \geq 0}\left\{ \frac{\lambda}{2} (r-R)^2 - g_f(R)\right\} r^{d-1}\rho(r) \dr \\ \leq & \int_{\R^d} S_\lambda \bar{f} (x) \mathrm{d} \mu (x).
	\end{align*}

\begin{theorem}\label{thm:radial barycenter}
	%	Let $\mu$ be the Wasserstein barycenter of $d$-dimensional radial contoured distributions $\{\mu_j = R_d(m_j, \rho_j)\}_{j=1}^N$ with weights $\{\lambda_j\}_{j=1}^N$, then $\mu$ is still a radial contoured distribution whose center is $m_* = \sum_{j=1}^N \lambda_j m_j$.
	Let $\{\mu_j = R_d(m_j,\rho_j)\}_{j=1}^N$ be $N$ radially contoured distributions and continuous generator functions $\{\rho_j\}_{j=1}^N$ are continuous. Then the Wasserstein barycenter $\mu$ of $\{\mu_j\}_{j=1}^N$ with weights $\lambda = (\lambda_1, \lambda_2, \cdots, \lambda_N) \in \Gamma_N$ is a generalized radial contoured distribution and the center is $m_* = \sum_{j=1}^N \lambda_j m_j$. 
\end{theorem}

\textbf{Proof.} \textbf{Step 1.} We first prove the case that $m_1 = m_2 = \cdots = m_N = 0$. Let $\mu$ be the barycenter of $\mu_j$ with weights $\{\lambda_j\}_{j=1}^N$. Suppose $\{f_j\}_{j=1}^N \subset X$ is a solution of (\ref{dual}), then $\sum_{j=1}^N f_j = 0$. It is easy to verify that $\bar{f_j} \in X$ and $\sum_{j=1}^N \bar{f_j}(x) = 0$. By Lemma \ref{lem3}, we have
	\begin{equation}
		\begin{aligned}
			F(f_1, f_2,\cdots, f_N) & = \sum_{j=1}^{N} \int_{\R^d} S_{\lambda_j} f_j(x) \mathrm{d} \mu_j(x) \\ & \leq \sum_{j=1}^N \int_{\R^d} S_{\lambda_j} \bar{f}(x) \mathrm{d} \mu_j(x) = F(\bar{f_1}, \bar{f_2}, \cdots, \bar{f_N}).
		\end{aligned}
	\end{equation}
	Therefore, $\{\bar{f_j}\}_{j=1}^N$ is also a solution of (\ref{dual}). Thus we can assume $f_j$ is itself radial. Then by Lemma \ref{lem2}, $S_{\lambda_j} f_j$ are radial functions. It follows by definition that $\phi_j$ (defined in (\ref{phi})) is also radial. Then by Theorem \ref{thm:barycenter characterization} and Lemma \ref{lem1}, $\mu = (\nabla \phi_1)_\# \mu_1$ is still a radial distribution.
	
	\textbf{Step 2.}	Let $\nu_j = R_d(0,\rho_j)$, and $\nu$ be their barycenter with weights $\{\lambda_j\}_{j=1}^N$. According to the previous proof, we have $\nu = R_d(0,\rho)$. Now we will prove $\mu = R_d(m_*,\rho)$. By Theorem \ref{thm:barycenter characterization}, there exist convex potentials $\{\psi_j\}_{j=1}^N$ such that $(\nabla \psi_j)_\# \nu_j = \nu$ and a constant $C$ such that
	\begin{equation}
		\sum_{j=1}^N \psi_j^*(y) \leq \frac{\norm{y}^2}{2} + C.
	\end{equation}
	Consider $\varphi_j(x) = \psi_j(x-m_j)+\inner{m_*}{x}$. Obviously $\varphi_j$ is convex. Since $(\nabla \psi_j)_\# \nu_j = \nu$, we have
	\begin{align*}
		\int_{\R^d} h(x) \mathrm{d} (\nabla \varphi_j)_\# \mu_j(x) & = \int_{\R^d} h(\nabla\psi_j(x-m_j)+m_*) \rho_j(\norm{x-m_j}) \dx \\
		& = \int_{\R^d} h(\nabla \psi_j(x) + m_*) \rho_j(\norm{x}) \dx = \int_{\R^d} h(x + m_*) \rho(\norm{x}) \dx \\
		& = \int_{\R^d} h(x) \rho(\norm{x-m_j} \dx) = \int_{\R^d} h(x) \mathrm{d} \mu(x), \quad \forall \, h \in C(\R^d),
	\end{align*}
	which means $\mu = (\nabla \varphi_j)_\# \mu_j$. Therefore $\nabla \varphi_j$ is the Monge map transporting $\mu_j$ to $\mu$ for every $j = 1, 2, \cdots, N$ by Theorem \ref{thm:Brenier}. Moreover
	\begin{align*}
		\varphi_j^*(y) & = \sup_{x \in \R^d} \left\{ \inner{x}{y} - \varphi_j(x) \right\} = \sup_{x \in \R^d} \{ \inner{x}{y} - \psi(x-m_j) - \inner{m_*}{x} \} \\
		& = \sup_{x \in \R^d} \left\{ \inner{x}{y-m_*} - \psi_j(x) \right\} + \inner{m_j}{y - m_*} = \psi_j^*(y-m_*) + \inner{m_j}{y-m_*}.
	\end{align*}
	Therefore
	\begin{align*}
		\sum_{j=1}^N \lambda_j \varphi_j^*(y) & = \sum_{j=1}^N \lambda_j \psi_j^*(y-m_*) +  \sum_{j=1}^N \lambda_j \inner{m_j}{y-m_*}\\
		& \leq \frac{\norm{y - m_*}^2}{2} + \inner{m_*}{y-m_*} + C \\ & = \frac{\norm{y}^2}{2} + C - \frac{\norm{m_*}^2}{2}.
	\end{align*}
	By Theorem \ref{thm:barycenter characterization}, $\mu$ is the barycenter of $\{\mu_j\}$ with weights $\{\lambda_j\}$.

\begin{remark}
	This theorem holds true when $\{\mu_j\}_{j=1}^N$ are $N$ generalized radially contoured distributions.
\end{remark}

When all the distributions $\{\mu_j\}_{j=1}^N$ share the same generator function, their Wasserstein barycenter also has a simple form.

\begin{example}\label{eg2}
	Let $\mu_j = R_d(m_j,c_j,\rho)$, $j = 1, 2, \cdots, N$, where $\rho$ is continuous. Then the barycenter of $\{\mu_j\}_{j=1}^N$ with weights $\lambda = (\lambda_1, \lambda_2, \cdots, \lambda_N) \in \Gamma_N$ is $\mu = R_d(m_*,c_*,\rho)$, where
	\begin{equation*}
		m_* = \sum_{j=1}^N \lambda_j m_j, \quad c_* = \sum_{j=1}^N \lambda_j c_j.
	\end{equation*}
\end{example}

\textbf{Proof.} According to Example \ref{eg1}, the potential $\psi_j$ such that $\nabla \psi_j$ is the Monge map transporting $\mu_j$ to $\mu$ is given by
	\begin{equation*}
		\psi_j(x) = \frac{1}{2} \frac{c_*}{c_j} \norm{x-m_j}^2 + \inner{m_*}{x}.
	\end{equation*}
	The Legendre conjugate of $\psi_j$ is
	\begin{equation*}
		\psi_j^*(y) = \sup_{x \in \R^d} \{\inner{y}{x} - \psi_j(x)\} = \frac{1}{2} \frac{c_j}{c_*}\norm{y}^2 + \inner{m_j - \frac{c_j}{c_*}m_*}{y} + C_j,
	\end{equation*}
	where $C_j$ is a constant independent of $y$. Then
	\begin{equation*}
		\sum_{j=1}^N \lambda_j \psi_j^*(y) = \frac{1}{2} \norm{y}^2 + \sum_{j=1}^N C_j.
	\end{equation*}
	By Theorem \ref{thm:barycenter characterization}, we get $\mu = R_d(m_*,c_*,\rho)$ is the barycenter of $\{\mu_j\}_{j=1}^N$ with weights $\{\lambda_j\}_{j=1}^N$.

\begin{figure}[htbp]
	\centering
	\includegraphics[width=0.3\textwidth]{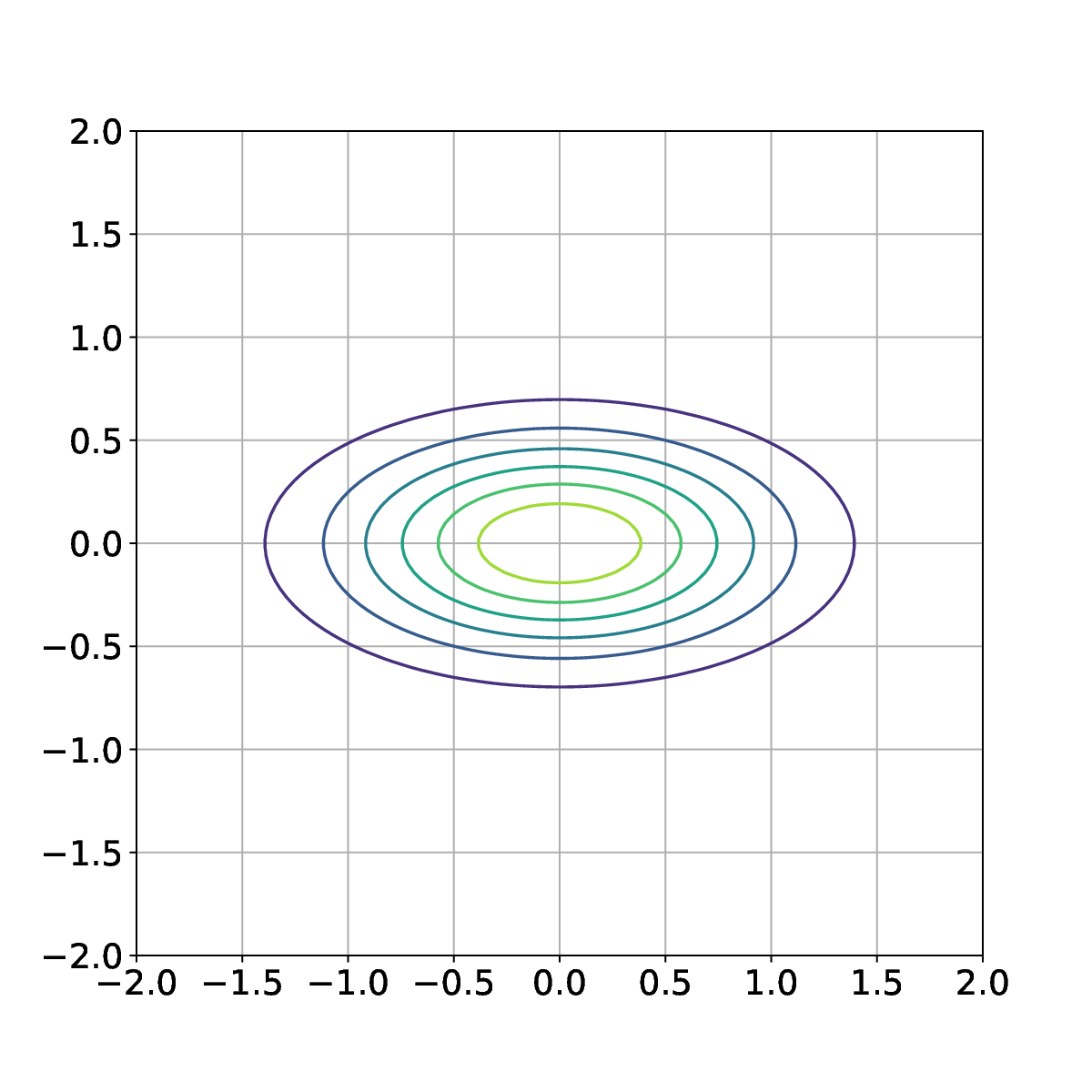}
	\includegraphics[width=0.3\textwidth]{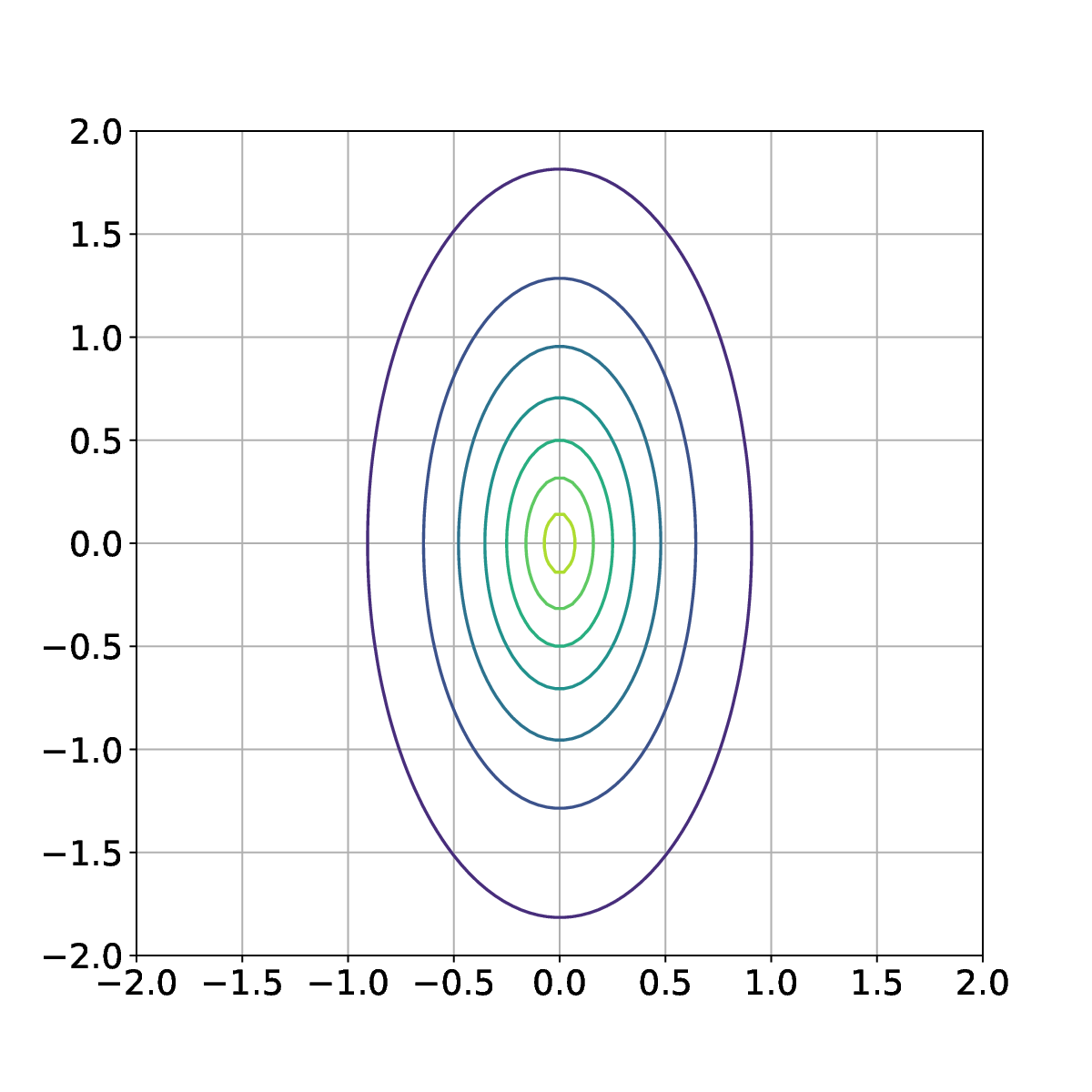}
	\includegraphics[width=0.3\textwidth]{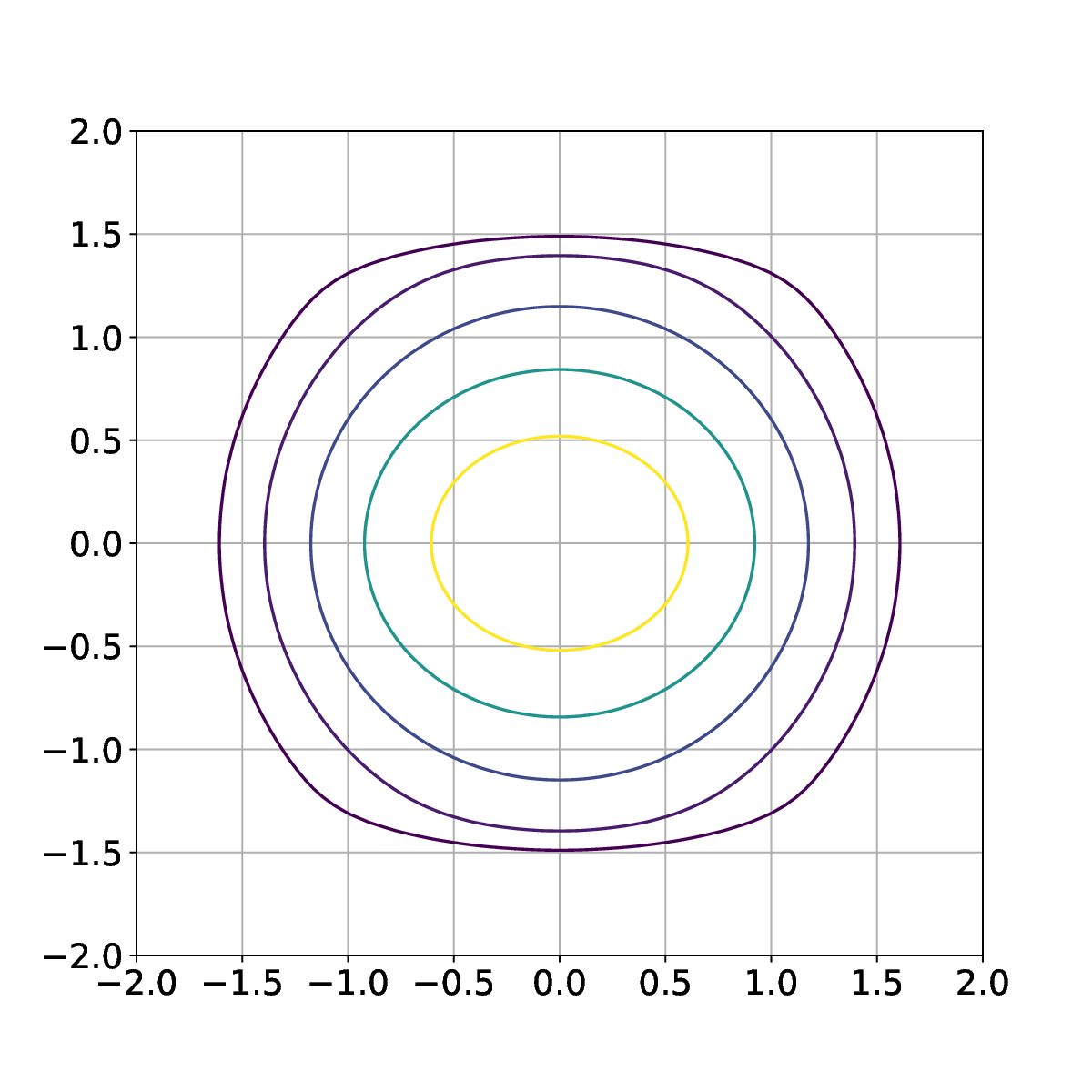}
	\includegraphics[width=0.3\textwidth]{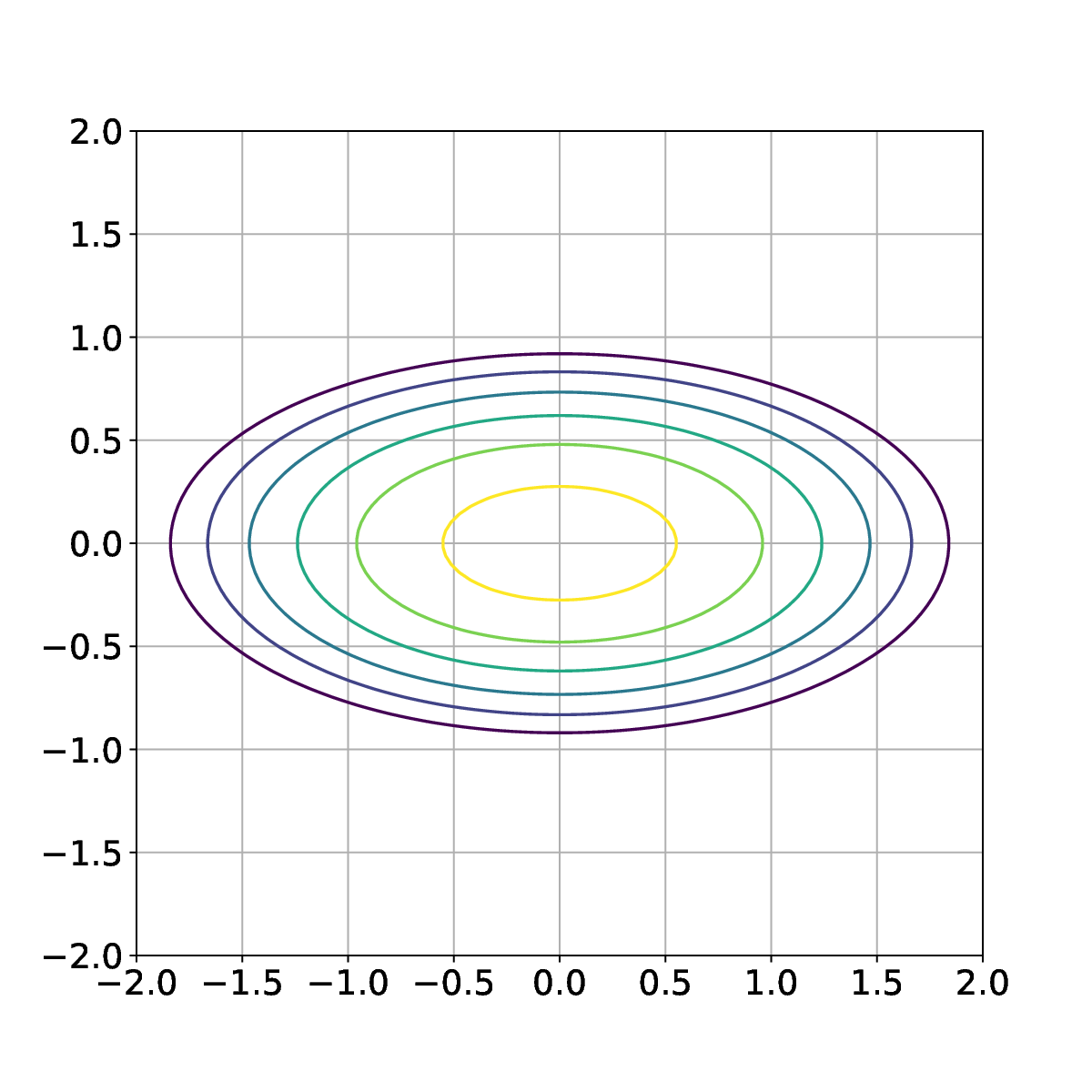}
	\includegraphics[width=0.3\textwidth]{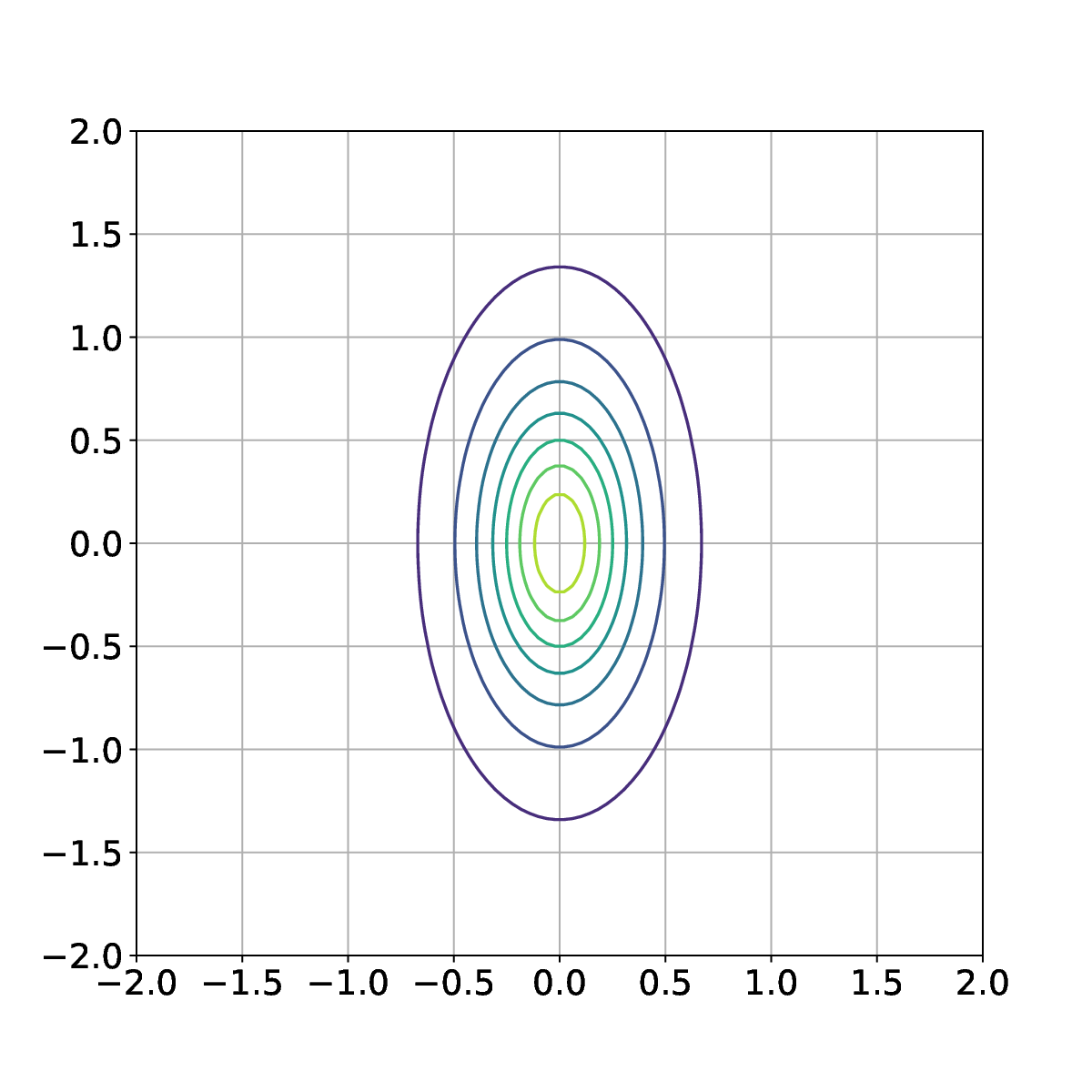}
	\includegraphics[width=0.3\textwidth]{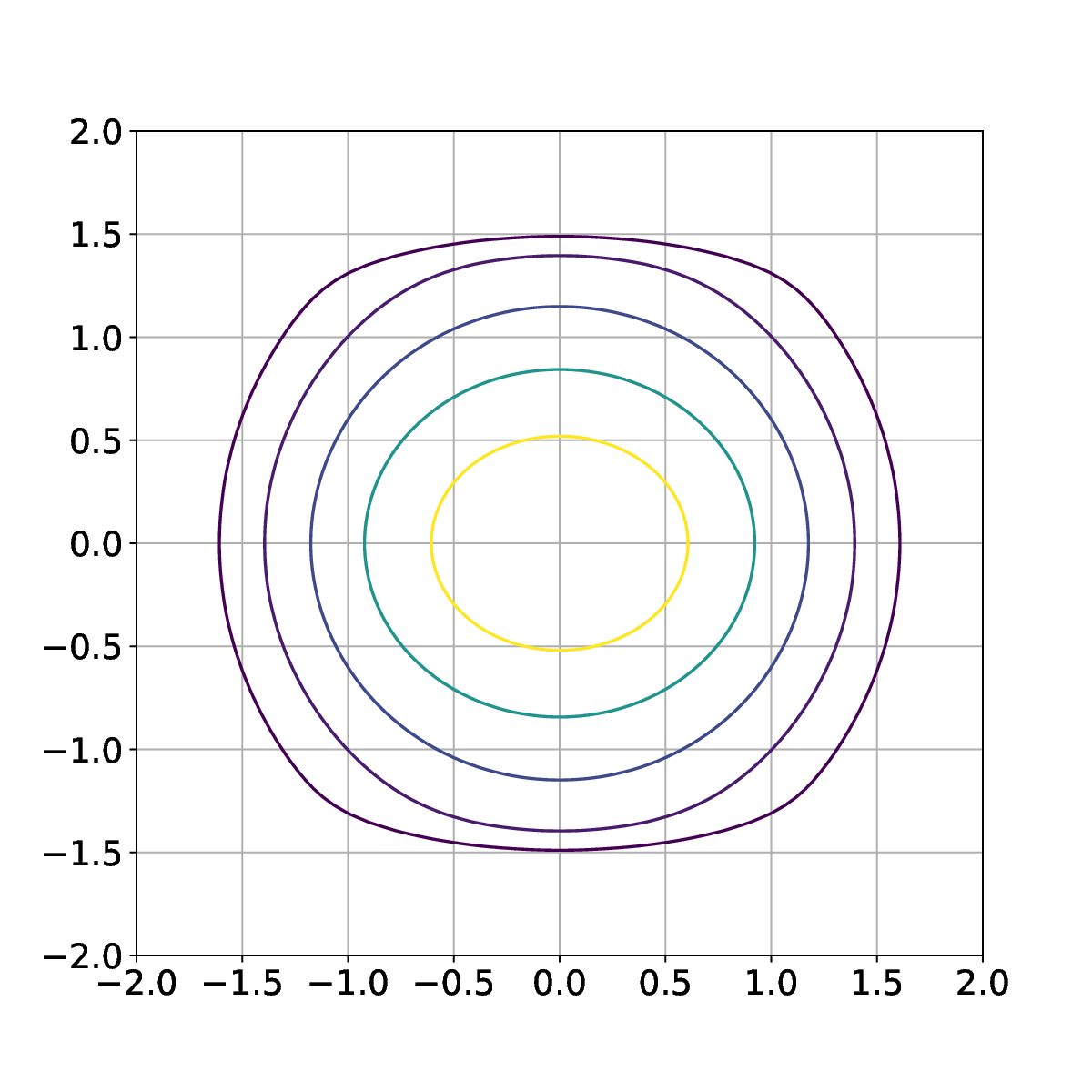}
	%	\subfigure{
		%		\includegraphics[width=0.3\textwidth]{ellip1.eps}
		%	}
	%	\subfigure{
		%		\includegraphics[width=0.3\textwidth]{ellip2.eps}
		%	}
	%	\subfigure{
		%		\includegraphics[width=0.3\textwidth]{bary.eps}
		%	} \\
	%	
	%	\subfigure{
		%		\includegraphics[width=0.3\textwidth]{ellip21.eps}
		%	}
	%	\subfigure{
		%		\includegraphics[width=0.3\textwidth]{ellip22.eps}
		%	}
	%	\subfigure{
		%		\includegraphics[width=0.3\textwidth]{bary2.eps}
		%	}
	
	\caption{The counterexamples that the Wasserstein barycenter of two elliptically contoured distributions may not be elliptically contoured. The first row is the first counterexample and the second row is the second one. In each row, the first two figures plot the contours of marginals, and the last figure shows the contours of the barycenter with weights $(0.5,0.5)$.}
	
	\label{fig counterexample}
\end{figure}

\section{Numerical counterexamples for elliptically contoured distributions}\label{sec5}

An elliptically contoured distribution $\mu$ on $\R^d$ takes the form
\begin{equation*}
	\mathrm{d} \mu(x) = \frac{1}{Z} \rho\left((x-m)^T \Sigma^{-1} (x-m)\right) \dx, 
\end{equation*}
where $\Sigma$ is a positive definite symmetric matrix, $\rho$ is a measurable function that maps $\R_+$ to $\R_+$ , and $Z$ is the normalizer. Such a distribution is denoted by $\mu = E_d(m,\Sigma,\rho)$. Similar to Theorem 6.1 in \cite{AguehMartial2011Bitw} for Gaussian distributions, the Wasserstein barycenter of $E_d(m_i,\Sigma_i,\rho)$ is of the form $E_d(m_*,\Sigma_*,\rho)$, which is still an elliptically contoured distributions. Although, we have proved that the Wasserstein barycenter for radially contoured distributions with different generator function retains radially contoured, similar result does not hold true for elliptically contoured distributions. We provide two numerical counterexamples here to show this conclusion in Figure \ref{fig counterexample}.

In the first example, the two elliptically contoured distributions are given by
\begin{equation*}\label{mu1}
	\mathrm{d} \mu_1 (x) = \frac{1}{Z_1} \exp\left(-(x_1^2 + 4 x_2^2)\right) 
	\mathrm{d} x_1 \mathrm{d} x_2,
\end{equation*}
\begin{equation*}\label{mu2}
	\mathrm{d} \mu_2 (x) = \frac{1}{Z_2} \exp\left(-(4 x_1^2 + x_2^2)^{0.6}\right) \mathrm{d} x_1 \mathrm{d} x_2.
\end{equation*}
In the second example, the two distributions are given by
\begin{equation*}\label{mu1h}
	\mathrm{d} \hat{\mu}_1 (x) = \frac{1}{\hat{Z}_1} \max\{(1-0.25x_1^2-x_2^2),0\}
	\mathrm{d} x_1 \mathrm{d} x_2,
\end{equation*}
\begin{equation*}\label{mu2h}
	\mathrm{d} \hat{\mu}_2 (x) = \frac{1}{\hat{Z}_2} (1+4x_1^2+x_2^2)^{-2} \mathrm{d} x_1 \mathrm{d} x_2.
\end{equation*}
%We solve the barycenter problem (\ref{Barycenter Problem}) for $\lambda_1 = \lambda_2 = 0.5$. The results are shown in Fig. \ref{fig counterexample}. In the first two figures, we plot the contours of the density functions of $\mu_0$ and $\mu_1$. In the right figure, we plot the contours of density function of the Wasserstein barycenter. We find that the inner contours seems like a circle or an ellipse. However, the outer two contours are far from ellipses.
We solve the barycenters problems numerically of the two set of elliptically contoured distributions with weights $\lambda_1=\lambda_2=0.5$. In the first row, we display the contours of the density functions of $\mu_1$, $\mu_2$ and their barycenter in the three figures from left to right. In the second row, the figures are contours of the density functions of $\hat{\mu}_1$, $\hat{\mu}_2$ and their barycenter. We find that for the contours of the barycenters, the inner contours seems like ellipses. However, the outer contours seems far from ellipses, especially the outermost one.

\section{Conclusion and discussion}
\label{sec6}

In this paper, we have studied the optimal transport and Wasserstein barycenter for radially contoured distributions with different generator functions. We derive the formulas of the Monge map, the McCann interpolation and Wasserstein distance of two radially contoured distributions. We prove that the Wasserstein barycenter of radially contoured distributions is still radially contoured. Furthermore, two numerical counterexamples are provided to show that the Wasserstein barycenter of elliptically contoured distributions is not necessarily elliptically contoured. 

Our Definition \ref{def Gradial} is not satisfactory. Here we allow the generator functions of radially contoured distributions to be measurable functions and a singular part in a generalized radially contoured distribution. But most of our results restrict the generator functions to be continuous. Meanwhile, when we use radial functions in applications, the radial functions are always continuous \cite{Light1992,Zhang1,Zhang2}. General measurable functions and Dirac-$\delta$ function will not be used. However, even if the generator functions are continuous, in our proof, we can not exclude the singular part in the barycenter. Whether the barycenter in this case is continuous or not is still a problem.

\section*{Acknowledgments}
This research is supported by National Key R\&D Program of China (2024YFA1012401), the Science and Technology Commission of Shanghai Municipality (23JC1400501), and Natural Science Foundation of China (12241103).

%Bibliography
\bibliographystyle{unsrt}  
\bibliography{references}

\begin{thebibliography}{10}

\bibitem{EvansL.C1999DEMf}
Lawrence~Craig Evans and Wilfrid Gangbo.
\newblock {\em {Differential Equations Methods for the Monge-Kantorovich Mass
  Transfer Problem}}, volume 137.
\newblock American Mathematical Society, Providence, 1999.

\bibitem{Villani2009OTOa}
Cédric Villani.
\newblock {\em Optimal Transport: Old and New}, volume 338 of {\em Grundlehren
  der mathematischen Wissenschaften}.
\newblock Springer-Verlag, Berlin, Heidelberg, 2009.

\bibitem{HakerSteven2004Omtf}
Steven Haker, Lei Zhu, Allen Tannenbaum, and Sigurd Angenent.
\newblock {Optimal Mass Transport for Registration and Warping}.
\newblock {\em International Journal of Computer Vision}, 60(3):225--240, 2004.

\bibitem{6502714}
Martin Mueller, Peter Karasev, Ivan Kolesov, and Allen Tannenbaum.
\newblock {Optical Flow Estimation for Flame Detection in Videos}.
\newblock {\em IEEE Transactions on Image Processing}, 22(7):2786--2797, 2013.

\bibitem{korotin2023neural}
Alexander Korotin, Daniil Selikhanovych, and Evgeny Burnaev.
\newblock {Neural Optimal Transport}.
\newblock In {\em ICLR 2023}, pages 1--34, 2023.

\bibitem{monge1781memoire}
Gaspard Monge.
\newblock {\em {M{\'e}moire sur la th{\'e}orie des d{\'e}blais et des
  remblais}}.
\newblock De l'Imprimerie Royale, Paris, 1781.

\bibitem{KantorovichL.V.2006Otto}
Leonid~Vitalyevich Kantorovich.
\newblock {On the Translocation of Masses}.
\newblock {\em {Journal of Mathematical Sciences (New York, N.Y.)}},
  133(4):1381--1382, 2006.

\bibitem{BrenierYann1991Pfam}
Yann Brenier.
\newblock {Polar Factorization and Monotone Rearrangement of Vector-Valued
  Functions}.
\newblock {\em Communications on Pure and Applied Mathematics}, 44(4):375--417,
  1991.

\bibitem{McCannRobertJ.1997ACPf}
Robert~J. McCann.
\newblock {A Convexity Principle for Interacting Gases}.
\newblock {\em {Advances in Mathematics}}, 128(1):153--179, 1997.

\bibitem{AguehMartial2011Bitw}
Martial Agueh and Guillaume Carlier.
\newblock {Barycenters in the Wasserstein space}.
\newblock {\em SIAM Journal on Mathematical Analysis}, 43(2):904--924, 2011.

\bibitem{SolomonJustin2015CWdE}
Justin Solomon, Fernando De~Goes, Gabriel Peyré, Marco Cuturi, Adrian
  Butscher, Andy Nguyen, Tao Du, and Leonidas Guibas.
\newblock Convolutional wasserstein distances: Efficient optimal transportation
  on geometric domains.
\newblock {\em ACM transactions on graphics}, 34(4):1--11, 2015.

\bibitem{PeyreGabriel2019Cot}
Gabriel Peyré and Macro Cuturi.
\newblock {Computational Optimal Transport}.
\newblock {\em Foundations and Trends in Machine Learning}, 11(5-6):1--257,
  2019.

\bibitem{GelbrichMatthias1990OaFf}
Matthias Gelbrich.
\newblock {On a Formula for the $L^2$ Wasserstein Metric between Measures on
  Euclidean and Hilbert Spaces}.
\newblock {\em {Mathematische Nachrichten}}, 147(1):185--203, 1990.

\bibitem{DeledalleCharlesAlban2018Idwg}
Charles~Alban Deledalle, Shibin Parameswaran, and Truong~Q. Nguyen.
\newblock {Image Denoising with Generalized Gaussian Mixture model Patch
  Priors}.
\newblock {\em {SIAM Journal on Imaging Sciences}}, 11(4):2568--2609, 2018.

\bibitem{GalerneBruno2017SOTi}
Bruno Galerne, Arthur Leclaire, and Julien Rabin.
\newblock Semi-discrete optimal transport in patch space for enriching gaussian
  textures.
\newblock In {\em Geometric Science of Information}, volume 10589 of {\em
  Lecture Notes in Computer Science}, pages 100--108. Springer International
  Publishing AG, Switzerland, 2017.

\bibitem{TeodoroAfonsoM.A.M.2015Sida}
Afonso~M.A.M. Teodoro, Mariana~S.C. Almeida, and Mário~A.T. Figueiredo.
\newblock Single-frame image denoising and inpainting using gaussian mixtures.
\newblock In {\em ICPRAM 2015: Proceedings of the International Conference on
  Pattern Recognition Applications and Methods}, volume~2, pages 283--288,
  2015.

\bibitem{DelonJulie2020Awdi}
Julie Delon and Agnès Desolneux.
\newblock {A Wasserstein-Type Distance in the Space of Gaussian Mixture
  Models}.
\newblock {\em SIAM Journal on Imaging Sciences}, 13(2):936--970, 2020.

\bibitem{ChenYukun2020AWDa}
Yukun Chen, Jianbo Ye, and Jia Li.
\newblock {Aggregated Wasserstein Distance and State Registration for Hidden
  Markov Models}.
\newblock {\em IEEE Transactions on Pattern Analysis and Machine Intelligence},
  42(9):2133--2147, 2020.

\bibitem{ChenYongxin2019OTfG}
Yongxin Chen, Tryphon~T. Georgiou, and Allen Tannenbaum.
\newblock {Optimal Transport for Gaussian Mixture Models}.
\newblock {\em {IEEE Access}}, 7:6269--6278, 2019.

\bibitem{dusson2023wassersteintypemetricgenericmixture}
Geneviève Dusson, Virginie Ehrlacher, and Nathalie Nouaime.
\newblock A wasserstein-type metric for generic mixture models, including
  location-scatter and group invariant measures, 2023.

\bibitem{CuestaalbertosJ.A.1993OCoM}
Juan~Cuesta Albertos, Ludger Ruschendorf, and Araceli~Tuero Diaz.
\newblock Optimal coupling of multivariate distributions and stochastic
  processes.
\newblock {\em Journal of Multivariate Analysis}, 46(2):335--361, 1993.

\bibitem{SantambrogioFilippo2015OTfA}
Filippo Santambrogio.
\newblock {\em Optimal Transport for Applied Mathematicians: Calculus of
  Variations, PDEs, and Modeling}, volume~87 of {\em Progress in Nonlinear
  Differential Equations and Their Applications}.
\newblock Springer Nature, Cham, 1st ed. 2015 edition, 2015.

\bibitem{GomezEusebio2003Asoc}
Eusebio Gómez, Miguel~A. Gómez-villegas, and J.~Miguel Marín.
\newblock {A Survey on Continuous Elliptical Vector Distributions}.
\newblock {\em {Revista Matemática Complutense}}, 16(1):345--361, 2003.

\bibitem{Light1992}
W.~A. Light.
\newblock {\em Some Aspects of Radial Basis Function Approximation}, pages
  163--190.
\newblock Springer Netherlands, Dordrecht, 1992.

\bibitem{Zhang1}
Yunxin Zhang.
\newblock {Reconstruct Multiscale Functions Using Different RBFs in Different
  Subdomains}.
\newblock {\em Applied Mathematics and Computation}, 189:893--901, 2007.

\bibitem{Zhang2}
Yunxin Zhang.
\newblock {Solving Partial Differential Equations by Meshless Methods Using
  Radial Basis Functions}.
\newblock {\em Applied Mathematics and Computation}, 185:614--627, 2007.

\end{thebibliography}

\end{document}